\documentclass[12pt]{iopart}
\usepackage{iopams}
\usepackage{latexsym}

\newtheorem{prop}{Proposition}      
\newtheorem{thm}{Theorem}
\newtheorem{lema}{Lemma}
\newtheorem{cor}{Corollary}
\newtheorem{defn}{Definition}

\newcommand{\epf}{\hfill   \mbox{$\Box $}}

\newcommand{\rf}{\mbox{${\mathbb R}$}}

\newcommand{\lel}[1]{\mbox{${\widetilde{#1}}$}}

\newcommand{\dd}{\mbox{${\rm d}$}}

\begin{document}

\title{On  the Time Average of the Autocorrelation Function in Hamiltonian Dynamics}
\date{November, 2008}

\author{Pavle Saksida and   Toma\v z Prosen  }
\address{Faculty of Mathematics and Physics,
University of Ljubljana, Jadranska 19, 
1000 Ljubljana, Slovenia}

\ead{Pavle.Saksida@fmf.uni-lj.si, Tomaz.Prosen@fmf.uni-lj.si}

\begin{abstract}
Rigorous lower bound on the time-average of the autocorrelation function of an arbitrary $L^1$ observable
is proven in terms of conserved quantities and ergodic decompositions of the Hamiltonian dynamics. Improvements
with respect to the bounds given by Mazur and Suzuki are discussed.
\end{abstract}


%


\maketitle

\section{Introduction}

Understanding the physical implications of ergodic and mixing properties of Hamiltonian dynamical 
systems is one of the main issues of non-equilibrium statistical mechanics. It is evident 
that a complete integrability in the sense of Liouville, or existence of even a small (non-complete) set 
of non-trivial global constants of the motion, is sufficient for a Hamiltonian system to be 
non-ergodic and non-mixing. It has been suggested by Mazur \cite{mazur}, and later followed up by 
Suzuki \cite{suzuki}, that 
the existence of conservation laws can also be connected to divergent transport coefficients expressed 
- due to a {\em linear response theory} of Green and Kubo (see, e.g. \cite{kubo,visscher}) - 
in terms of the integrated time auto-correlation of the current observable.
Namely, Green-Kubo formula expresses the conductivity $\kappa_A$ with respect to a certain current observable 
$A$ (say, the energy-current for the thermal conductivity) as
$$
\kappa_A = \beta \int _0 ^\infty \langle A^0, A^t\rangle _{\beta} \ 
\dd t ,
$$
where $\beta$ is the inverse temperature, $A^t$ denotes the time evolution of $A$.
More precise meaning of the
notation will be explained later.

Mazur has pointed out an inequality between the time-integrated or time-averaged
auto-correlation of an arbitrary observable 
(which may also be interpreted thermodynamically as an {\em isothermal susceptibility}), 
and an algebraic expression which depends solely on the overlaps between the observable in question and 
the conservation laws expressed in terms of phase space integrals. His theorem is essentially a straightforward
consequence of an ergodic theorem of Khinchine \cite{khinchine} and demonstrates that completely integrable 
systems should generically behave  as ideal (ballistic) conductors of heat, electricity, etc.

Nevertheless, as such a result might seem quite natural from the point of view of dynamical systems and 
ergodic theory, it raised a lot of surprise and attention in the solid state community \cite{zotos,zotosreview}. Namely, 
as expressed in the language of  solid state physics, one-dimensional completely integrable many-particle systems 
(with a large number of, or in the thermodynamic limit, an infinite number of degrees of freedom) generically posses a finite Drude weight, thus a divergent zero-frequency (d.c.) conductivity, and therefore behave as 
{\em ballistic conductors}, at {\em all temperatures}.
Ballistic transport in a strongly interacting system at a non-vanishing temperature is certainly a statement 
which raises eyebrows of a solid state physicist.
However,   complete integrability is not even necessary to have such a striking physical 
implication, it is enough that at least one non-trivial conservation law exists which overlaps with 
the current \cite{prosen98}.

In this paper we present some results which provide sharp and easily computable bounds
of the autocorrelation function of an observable in terms of the conserved quantities of the system. 
Mazur's theorem can be considered as a linear (or first order) case of our bounds.
Under certain conditions, which are quite likely to be fulfilled in many concrete cases,
our results actually yield precise values of the time average of the autocorrelation function. 
Mazur's proof relies heavily on a rather involved probabilistic and measure-theoretic 
statement of Khinchine \cite{khinchine}. Our method of proof is completely different,
rather simpler and it avoids deep results of Khinchine entirely. Nevertheless, this method
seems to give  a good handle on the problem. In particular, it explains how  the ergodic properties
of the system on the one hand, and the form of the observable in question on the other,
affect the sharpness of the bounds. We expect that  saturability
of our bounds can be an interesting test (probe) of complete integrability. We intend to address
this issue in another paper.

In order to have geometric concepts well defined,  we shall assume that the number
of degrees of freedom is finite. Our results are essentially valid for the systems with the
infinite number of degrees of freedom,  which could be argued by taking the thermodynamic limit.

Our setup is  the following.  Let $(M, \omega, H)$ be a Hamiltonian system,
where $M$ is the phase space with the symplectic structure  $\omega$ and the
Hamiltonian $H$. Let 
\[
A : M \longrightarrow \rf
\]
be an arbitrary observable. By $A^t$ we denote its time evolution given by
\[
A^t(m) =  A(\gamma _m(t)) : M \longrightarrow \rf,
\]
where $\gamma _m(t) \colon \rf \to M$ is the solution of the system $(M, \omega, H)$
satisfying the initial condition $\gamma _m(0) = m$.

The object of our interest in this paper will be the time average 
of the autocorrelation function
\[
{\cal C}(A) = \lim _{T \to \infty} \frac{1}{T}
\int _0 ^T \langle A^0, A^t\rangle _{\beta} \ 
\dd t .
\]
The symbol  $\langle - , - \rangle _\beta$  denotes
the inner product on the $L^2$- space $L^ 2_ \beta(M)$ given by
\[
\langle F, G \rangle _\beta =  \frac{1}{{\cal Z}(\beta)} \int _M F(m) \cdot G(m) \ e^{- \beta H(m)} \ \dd m,
\]
and
\[
{\cal Z}(\beta) = \int _M e^{- \beta H(m)} \ dm
 \]
 is the partition function which we shall assume to exist for strictly positive values of
 $\beta$ (that is, for all finite temperatures if the system belongs to the realm of  statistical physics). 
 The measure $\dd m $ is given by
 the top exterior power $\omega^r$ of the symplectic form, where $r$ is the number of degrees
 of freedom.
In local canonical coordinates:
\[
\dd m = \dd q_1 \cdots \dd q_r  \dd p_1 \cdots  \dd p_r.
 \]

The purpose of this paper  is to give simple geometric  descriptions of ${\cal C}(A)$ 
in terms of a set of, say $k$, functionally independent conserved quantities
\[
H = H_1, H_2, \ldots, H_k : M \longrightarrow \rf.
\]
Our main result is the following bound on the  time average of the autocorrelation function:
\begin{equation} \fl
{\cal C}(A)  \geq \sum _{ {\bf l } = (l_1, \ldots ,l_k),  \  0 \leq |{\bf l}| \leq d  \atop   
{\bf n} = (n_1, \ldots , n_k),  \  0 \leq |{\bf n}|  \leq d} 
\langle A,  H_1^{l_1} \ldots H_k^{l_k} \rangle _\beta \ ({\bf H}^{-1})_{\bf{l}, \bf{n}} \
\langle A,  H_1^{n_1} \ldots H_k^{n_k} \rangle _\beta,
\label{new}
\end{equation}
where $|{\bf n}| = \sum _{i= 1}^{k} n_i$.
The elements of the matrix ${\bf H}$ are given by inner products
\[
{\bf H}_{{\bf l }, {\bf n }} = \langle H_1^{l_1} H_2^{l_2} \ldots H_k^{l_k}  \ ,  \ 
 H_1^{n_1} H_2^{n_2} \ldots H_k^{n_k} \rangle _\beta 
 \]
 in a suitably chosen order. Immediately after proving the above bound, we identify
 the conditions under which the bound becomes an equality. We note  that the order 
 $d$ can be infinite.
 
Note that the strict positivity of ${\cal C}(A)$, implied by the strict positivity of 
the right-hand-side of (\ref{new}),
implies ballistic transport, if $A$ is the corresponding current observable in the Green-Kubo theory. 
 
 The inequality (\ref{new})  is a non-linear  improvement of the inequality  
 \begin{equation}
 {\cal C}(A) = \lim _{T \to \infty} \frac{1}{T} \int _0 ^T \langle A^0, A^t \rangle _\beta \  \dd t
 \ \geq \ \sum _{i , j = 1}^k \langle A, H_i \rangle _\beta \cdot ({\bf H}^{-1})_{i, j} \cdot
 \langle A, H_j\rangle _\beta ,
 \label{Mazur1}
 \end{equation}
 originally given by P. Mazur in \cite{mazur}.
 Above ${\bf H}$ denotes the matrix with entries 
 \begin{equation}
 {\bf H}_{i, j} =  \langle H_i, H_j \rangle _\beta
 \label{Hmatrix}
 \end{equation}
 and ${\bf H}^{-1}$ is its inverse.
 
 In \cite{mazur}, working in the framework of statistical
 mechanics, Mazur treats the evolution $t \mapsto A(\gamma _m(t)) = A^t$ as 
 a stochastic process. The main tool he uses is the power spectrum $I(\omega )$
 of the process $A^t$. The power spectrum is given by the Fourier transform
 of the correlation function
 \[
 \phi (t) = \langle A^0, A^t \rangle _\beta,
 \]
 that is by
 \[
 I(\omega ) = \frac{1}{2 \pi} \int _{- \infty}^{\infty} \phi (t) \ e^{- i t \omega} \ \dd t . 
 \]
 The  essential ingredient of  Mazur's proof is the well-known
 result of Khinchine (\cite{khinchine}) which states that for every
 function $B \colon M \to \rf$ we have
 \[
 {\cal C}(B) = \lim _{T \to \infty} \frac{1}{T} \int _0 ^T \langle A^0, A^t \rangle _\beta =
 I(0 + ) - I (0 - ) \geq 0 .
 \]
 From the above, Mazur's inequality follows almost immediately. The analytical crux of
 the matter is indeed contained in the above deep result of Khinchine.
 
 Our treatment relies on direct and simple geometric considerations and does not resort to the 
 theory of stochastic processes, let alone to the Khinchine's result. Despite its simplicity, our
 geometric approach enables us to get a better grip on the internality. We give  simple 
 and meaningful conditions for the inequality to be saturated. We also show how the ergodic
 properties of the system $(M, \omega, H)$ affect the quantity ${\cal C}(A)$.

Throughout  this paper the observables $A \colon M \to \rf$ will be assumed to be elements
of $L^1_\beta(M)$,  where in general 
$L^p_\beta(M) $ is the $L^p$-space generated by the measurable real functions
$f \colon  M \to \rf$ which satisfy the integrability condition
\[
\int _M | f | ^p \ e^{-\beta H(m)} \dd m \ < \infty.
\]
Assuming that the partition function ${\cal Z}(\beta)$ is defined, the measure
$e^{-\beta H(m)} \dd m$ is finite. Therefore, we have the inclusion
$L_\beta ^1(M) \subset L_\beta^2(M)$, and thus $A$ will also be an element
of $L_\beta^2(M)$.

We will divide the discussion into two parts. Firstly, in section 2 we shall study
the so-called  ergodically regular case. This means that the dynamics is ergodic on the
joint level sets of the conserved quantities.  There the bound (\ref{new}) is proven
in theorem \ref{theorem2}, whereas in theorem \ref{theorem1} we 
give a useful expression of ${\cal C}(A)$ as the $L_2$-norm of a certain projection
of the observable $A$.  Theorem \ref{corsat} identifies the situation in which the 
bound (\ref{new}) is saturated.  It is reasonable to expect that the result of the
theorem \ref{corsat} will be useful in the context of the algebraically integrable
systems, where the conserved quantities are expressible as polynomials or analytic 
functions. We list some sufficient conditions for the saturation in the remarks
following theorem \ref{corsat}.
Secondly, in section 3 we study the general case, without the assumption 
of  ergodic regularity. The central technical result here is lemma \ref{lemma}, which
is needed to prove the bound (\ref{new}), (see corollary \ref{corlast}).  In
theorem \ref{theorem3} we prove that ${\cal C}(A)$ is equal to the 
$L_2$-norm of the orbital average of the observable $A$.
In proposition \ref{prop4} we show how ergodic decompositions can be used 
to further improve the bound (\ref{new}).

\section{The  ergodically regular case}

Let, as above, $(M, \omega , H)$ be a Hamiltonian system with $k < 2 r =  \ {\rm dim} (M)$
{\it functionally independent}
conserved quantities
\[
H = H_1, H_2, \ldots, H_k : M \longrightarrow \rf
\]
which are not necessarily in involution. Let 
\begin{eqnarray}
{\cal H}   =   (H_1, \ldots , H_k)  :  & \ M \longrightarrow D \subset \rf ^k \\ \nonumber 
& \  m \longmapsto   {\cal H}(m) = ( \alpha _1, \ldots, \alpha _k)
\end{eqnarray}
denote the moment map.  We shall first consider the systems whose ergodic
behaviour is simple in the sense that it is completely determined by the conserved quantities. 

For the sake of brevity we shall denote the level sets of the
moment map by 
\[
{\cal L}_\alpha = {\cal H}^{-1}(\alpha).
\]
\begin{defn}
The system  $(M, \omega, H)$ is called  ergodically regular,  if it is ergodic on level sets
${\cal L}_\alpha$ for almost every $\alpha \in D$. More precisely, for almost every
$\alpha \in D$ and for almost every $m \in {\cal L}_\alpha$ we have
\[
\lim _{T \to \infty } \frac{1}{T} \int _0 ^T B(\gamma _{m} (t)) \ \dd t =
\frac{1}{{\rm Vol}({\cal L}_\alpha)} \int _{{\cal L}_\alpha} B(m_\alpha) \ \dd m_\alpha,
\]
where $B \in L^1_\beta(M)$ is an arbitrary observable, and $\dd m_\alpha$
is the generalized microcanonical measure on ${\cal L}_\alpha $ induced 
by the  measure
$\dd m$  which is given by $\omega ^r$.
\end{defn}
A general example of  ergodically regular systems are Liouville integrable systems.
For the definition of the induced microcanonical measure, see \cite{minlos}.
In the above definition it is assumed that  every volume ${\rm Vol}({\cal L}_\alpha)$ is
finite, but a moment of  thought shows that this is a consequence 
of the finiteness of the partition function ${\cal Z}(\beta)$.

The key ingredient in the study of  ergodically regular systems is the averaging map 
\begin{equation}
B   \  \stackrel{ {\cal T}}{\longmapsto}  \  B^{\cal H} 
= \frac{1}{{\rm Vol}({\cal L}_\alpha)} \int _{{\cal L}_\alpha} B(h_\alpha) \ \dd h_\alpha .
\label{trans}
\end{equation}
Let
$L^1_\beta(D)$ be the $L^1$-space of measurable functions on $D$ with respect to the 
measure $\mu_D$ given by
\begin{equation}
\dd \mu _D =   \frac{1}{{\cal Z}(\beta)} \
{\rm Vol}({\cal L}_\alpha )\ e^{-\beta H^{\cal H}(\alpha)} \ \dd \alpha = 
\frac{1}{{\cal Z}(\beta)} \ {\rm Vol}({\cal L}_\alpha )\ e^{-\beta \alpha _1} \ \dd \alpha .
\label{Dmeas}
\end{equation}
Then the  assignment (\ref{trans}) 
defines the operator
\begin{equation}
\begin{array}{cccc}
{\cal T} : & L^1_\beta(M)  & \longrightarrow &  L^1_\beta(D)   \\
       & & &   \\
                 & B & \longmapsto & B^{\cal H}.
  \end{array}
  \end{equation}
 Since our measures on $M$ and $D$ are both finite, we have $L_\beta ^2(M) \subset L_\beta ^1(M)$
 and $L_\beta^2(D) \subset L_\beta ^1(D)$.  Obviously, the averaging operator descends to the map
 \[
 {\cal T} : L_\beta ^2(M)  \longrightarrow L_\beta^2(D)
 \]
 between the $L^2$-spaces.
 Let now the operator
 \[
 {\cal R} : L_\beta^i(D) \longrightarrow L_\beta^i(M), \quad i = 1,2
 \]
 be given by
 \[
 ({\cal R}(G))(m) = G({\cal H}(m)).
 \]
 The composed operator
 \begin{equation}
 {\cal P} = {\cal R} \circ {\cal T} : L_\beta ^2(M)  \longrightarrow L_\beta ^2 (M)
 \label{proj}
 \end{equation}
 is then clearly a projector.
 \begin{prop}
 The operator ${\cal P}$ given by (\ref{proj}) is a continuous orthogonal operator 
 and therefore  symmetric  on $L_\beta^2(M)$.
 \end{prop}
 {\bf Proof.} To prove the orthogonality, let $F$ be an element of the kernel of ${\cal P}$, 
 and let $G$ lie in its image. Then
 
 \begin{eqnarray*} \fl
 \int _M F(m) \ G(m) \ e^{-\beta H(m)} \ \dd m  & = &
 \int _D \ \dd \alpha \int _{{\cal L}_\alpha} F(m_\alpha) \ G(m_\alpha ) \ 
 e^{-\beta H(m_\alpha)} \ \dd m_\alpha \\ \nonumber
 & = & \int _D G^{\cal H}(\alpha) \ e^{-\beta H^{\cal H}(\alpha)} \Bigl(
 \int _{{\cal L}_\alpha} F(m_\alpha) \ \dd m_\alpha \Bigr) \ \dd \alpha \\
 & = & 0 .
 \end{eqnarray*}
 Considering suitable convergent sequences in $L_\beta ^2(M)$ it is easily seen that the subspaces
 ${\rm Ker}({\cal P})$ and ${\rm Im}({\cal P})$ are closed, which implies that ${\cal P}$ is a continuous
 orthogonal projection and therefore a symmetric map of $L_\beta^2(M)$ into itself.
 
 \epf
  
  From the above proposition we immediately obtain the basic geometric description 
  of the time average ${\cal C}(A)$ in the  ergodically regular case.
  By $\| - \|_D$ we shall denote the norm on $L_\beta^2(D)$ with respect to the measure
  $\dd \mu _D$
  introduced above. 
  
  \begin{thm}
  Let $A$ be an arbitrary measurable and integrable observable of the Hamiltonian system
  $(M, \omega, H)$. Then for the time average of its correlation function we have
  \[
  {\cal C}(A) =  
  \|{\cal P}(A) \|^2_\beta = \| A^{\cal H} \|^2_D,
  \]
  where $\| - \|_D$ denotes the $L_2$-norm on $L_\beta^2(D)$ given by
  the measure (\ref{Dmeas}).
  \label{theorem1}
  \end{thm}
{\bf Proof.}
The fact that our system is  ergodically regular 
gives
\begin{eqnarray*} \fl
    \lim _{T \to \infty}
\frac{1}{T}  \ \int _0^T  \langle A^0, A^t \rangle _{\beta} \ \dd t  &  =  &  
   \lim _{T \to \infty} \frac{1}{T}  \ \int _0^T  \frac{1}{{\cal Z}(\beta)} \
\int _M  A(\gamma _m(0)) \cdot  A(\gamma _m(t)) 
\cdot e^{- \beta H(m)}  \ \dd m                   \\
   & =  &
\frac{1}{  {\cal Z}(\beta)} \int _M \Bigl( A(\gamma _m (0)) \cdot 
\lim _{T \to \infty} \frac{1}{T}
\int _0^T A(\gamma _m(t)) \  \cdot e^{- \beta H(m)} \ dt \Bigr) \  \dd m   \\
   & = &
\frac{1}{  {\cal Z}(\beta)} \int _M \Bigl( A(m) \cdot
\frac{1}{{\rm Vol({\cal L}_{{\cal H}(m)}})} \ \int _{{\cal L}_{{\cal H}(m)}} 
A(h)   \ \dd h \ \Bigr)
\ e^{- \beta H(m)}\ \dd m \  \\
& =  & \frac{1}{  {\cal Z}(\beta)} \int _M  A(m) \ {\cal P}(A)(m) \
e^{- \beta H(m)}  \ \dd m .
\end{eqnarray*}
Using our notation and the fact that ${\cal P}$ is a symmetric projection
we  get the first equality of the theorem,
\[
\langle A, {\cal P}(A)\rangle _\beta = \langle A, {\cal P}^2 (A) \rangle _\beta =
\langle {\cal P}(A), {\cal P}(A) \rangle _\beta = \|{\cal P}(A)\|^2_\beta.
\]
Finally, for every pair of functions $A, B \colon M \to \rf$ we have
\begin{eqnarray} \fl
\langle {\cal P}(A), {\cal P}(B) \rangle _\beta & = &
 \frac{1}{{\cal Z}(\beta)} \  \int _M  {\cal P}(A)(m)\  {\cal P}(B)(m)
\ e^{-\beta H(m)} \ \dd m  \nonumber \\
& =  & \frac{1}{{\cal Z}(\beta)} \ \int _D \dd \alpha 
\int _{{\cal L}_\alpha} {\cal P}(A)(m_\alpha)\  {\cal P}(B)(m_\alpha)
\ e^{-\beta H(m_\alpha)} \ \dd m_\alpha \\
& =  &\frac{1}{{\cal Z}(\beta)} \ 
\int _D A^{\cal H}(\alpha ) \ B^{\cal H}(\alpha)  \ e^{-\beta \alpha _1}
{\rm Vol}({\cal L _\alpha}) \ \dd \alpha    \nonumber \\ 
& = & 
\langle A^{\cal H}, B^{\cal H} \rangle _D, \nonumber
\label{dbeta}
\end {eqnarray}
which proves our second equality.

\epf

We shall now turn to the analogues of the right-hand side of Mazur's inequality. These
expressions will give us estimates for ${\cal C}(A)$. In concrete cases these expressions
are likely to be simpler to calculate than the norm of $A^{\cal H}$. More  importantly,
they will be useful in the treatment of the general, ergodically irregular case.

Let us consider the set of all monomials  corresponding to the multi-indices
${\bf n} = (n_1, \ldots, n_k)$ of degree $d$ or less
\[
H_1^{n_1} H_2^{n_2} \ldots H_k^{n_k}, \quad \   0 \leq |{\bf n}| = \ n_1 + n_2 + \ldots n_k  \leq d
\]
composed of the conserved quantities
$H = H_1, H_2, \ldots , H_k$ of our system.  We
shall order these monomials by a combination of the ordering
by degree and the lexicographical ordering,
\begin{equation} \fl
i\!f \ |{\bf n} _1| < |{\bf n }_2|, \ then \  {\bf n}_1 < {\bf n}_2\ ; \quad \
i\!f \ |{\bf n} _1| = |{\bf n }_2|, \ then \ ordered \ lexicographically.
  \label{ord}
  \end{equation}
By $V_d$ we shall denote the subspace of $L_\beta^2(M)$ spanned by the 
above monomials,
\[
V_d = {\rm span} \{H_1^{n_1} H_2^{n_2} \ldots H_k^{n_k}; \  |{\bf n}| \leq d\}
\subset L_\beta^2(M) .
\]
The fact that the conserved quantities are functionally independent implies the 
linear independence of the  monomials.
The basis $\{H_1^{n_1} \ldots H_k^{n_k}\}$   of $V_d$ is not orthonormal, therefore 
we will need the matrix of the inner product $\langle - , - \rangle _\beta $ on $V_d$
corresponding to our basis.  The elements of the  inner product matrix ${\bf H}$
are given by
\begin{eqnarray*}
({\bf H}) _{o({\bf l}), o({\bf n})} & = & \langle H_1^{l_1} H_2^{l_2} \ldots H_k^{l_k} \ , 
\ H_1^{n_1} H_2^{n_2} \ldots H_k^{n_k} \rangle _\beta  \\
& = & \int _M H_1^{l_1 + n_1}(m) H_2^{l_2 + n_2}(m) \ldots H_k^{l_k + n_k}(m) \ 
e^{-\beta H(m)} \ \dd m,
\end{eqnarray*}
where $o({\bf l })$, $o({\bf n})$ are integers given by the ordering (\ref{ord}).
The matrix ${\bf H}$ is non-singular due to the linear independence of the monomials.
\begin{thm}
Let $A$ be an observable on an ergodically regular Hamiltonian system $(M, \omega, H)$
with the additional conserved quantities $H_2, \ldots H_k$. Then
for every positive integer $d$ we have
\begin{equation} \fl
{\cal C}(A)   \geq \sum _{o({\bf l}), o({\bf n }) = 0}^{\nu (d)}
\langle A  ,   H_1^{l_1} H_2^{l_2} \ldots H_k^{l_k} \rangle _\beta \ 
({\bf H} ^{-1})_{o({\bf l}), o({\bf n})} \ 
\langle A  ,   H_1^{n_1} H_2^{n_2} \ldots H_k^{n_k} \rangle _\beta,
 \label{new1}
 \end{equation}
 where $\nu (d)$ denotes the number of different monomials of degrees ranging
 between $0$ and $d$ in $k$ unknowns.
 \label{theorem2}
  \end{thm}
  The bound (\ref{new1}) could be of practical importance, since it is relatively easily calculable
  in many cases. We notice that the left-hand side of (\ref{new}) is independent of $d$, therefore
  theorem \ref{theorem2}  has the following immediate corollary.
  \begin{cor}
  Let $(M, \omega, H)$ and $A \colon M \to \rf$ be as above. Then
  \begin{equation} \fl
  {\cal C}(A)   \geq \sum _{o({\bf l}), o({\bf n }) = 0}^{\infty}
  \langle A  ,   H_1^{l_1} H_2^{l_2} \ldots H_k^{l_k} \rangle _\beta \ 
  ({\bf H} ^{-1})_{o({\bf l}), o({\bf n})} \ 
  \langle A  ,   H_1^{n_1} H_2^{n_2} \ldots H_k^{n_k} \rangle _\beta.
   \label{infinite}
 \end{equation}
  \end{cor}
  \epf
  
  \noindent 
  Clearly, the bound (\ref{infinite}) is sharp. From the practical point of view, it is in general less
  useful. The main problem is the evaluation of the inverse of the infinite matrix
  ${\bf H}$. This can be tackled by replacing the monomials by products of
  suitably scaled orthogonal polynomials.

  \noindent {\bf Proof of Theorem \ref{theorem2}.} Let 
  ${\cal U}= {\cal P}(L_\beta ^2(M)) \subset L_\beta ^2(M)$
  be the image of the projector ${\cal P}$.  This is a closed subspace
  in $L_\beta ^2(M)$ and the inherited inner product gives it the structure of
  a Hilbert space. Let
  \[ 
  \pi :  {\cal U} \longrightarrow V_d
  \]
  be the orthogonal projection. We shall prove that the expression on the right-hand
  side of (\ref{new1}) is equal to the norm of the vector $\pi ({\cal P}(A))$.
  Since by Pythagoras' theorem any orthogonal projection of a vector is
  shorter than the vector itself, (\ref{new1}) will   follow immediately from
  theorem \ref{theorem1}.
  
  Let us introduce a shorter notation
  \[
 h_{o({\bf n})} =   H_1^{n_1} H_2^{n_2} \ldots H_k^{n_k}  
   \]
  for the elements of the basis of $V_d$, and let $\{h^*_i;  i = \nu(d) + 1,  \nu(d) + 2,\ldots
  \} $ 
  be the basis of the orthogonal complement $V_d^{\perp} \subset {\cal U}$ 
  composed of the vectors from the basis dual to $\{h_n\}_{n\in \mathbb{N}}$. This means
  \[
  \langle h_i^* , h_j \rangle _\beta = 0, \quad for \quad i = 0 \ldots , \nu (d), \quad j = \nu (d) + 1, 
  \nu (d) + 2, \ldots  .
  \]
  Let
  \begin{equation}
  {\cal P}(A) = \sum _{i = 0}^{\nu(d)} c_i  \ h_i \ + \ \sum _{j = \nu(d) + 1}^{\infty} d_j  \  h_j^*
  \label{decomp}
  \end{equation}
  be the orthogonal decomposition of ${\cal P}(A)$ with respect to ${\cal U} =
  V_d \oplus V_d^{\perp}$. Taking the inner product of this expression with each of
  $h_k$ for $\ k = 1, \ldots , \nu(d)$ gives
  \[
  \langle {\cal P}(A), h_k \rangle _\beta =  \sum _{i = 0}^{\nu(d)} c_i  
  \langle  h_i, h_k \rangle _\beta.
  \]
  Inverting the matrix ${\bf H}_{i, k} =  \langle h_i, h_k \rangle _\beta$ yields the expression
  of the coefficients $c_i$:
  \begin{equation}
  c_i = \sum _{j = 0}^{\nu(d)} ({\bf H}^{-1})_{i, j} \ \langle {\cal P}(A), h_j \rangle _\beta .
  \label{coeffs}
  \end{equation}
  The inner product matrix of $V_d \subset L_\beta^2(M)$ with respect to the basis
  $\{h_i, \ldots , h_{\nu (d)}\}$ is $\bf{H}$, thus
  the squared norm $\| \pi ({\cal P}(A)\|^2_\beta$ of the orthogonal projection of ${\cal P}(A)$
  onto $V_d$ is equal to
  \[
\| \pi ({\cal P}(A)\|^2_\beta = (c_0, \ldots , c_{\nu (d)}) \cdot {\bf H} 
\cdot \pmatrix{ c_0 \cr
\vdots \cr
c_{\nu (d)} \cr} . 
\]
 By means of (\ref{coeffs}) we can express this quantity in terms of the products
 $\langle {\cal P}(A), h_i \rangle _\beta$. The inner product matrix ${\bf H}$ 
 is replaced by ${\bf H}^{-1} \cdot {\bf H} \cdot {\bf H}^{-1}$, so we
 get
\[ \fl
\| \pi({\cal P}(A))\|^2 
 = {\bf (}\langle {\cal P}(A), h_0\rangle _\beta ,  
 \ldots ,  \langle {\cal P}(A), h_{\nu (d)}\rangle _\beta  {\bf )} \cdot
 {\bf H}^{-1}  \cdot
    \pmatrix{\langle {\cal P}(A),  h_0\rangle _\beta \cr
                                            \vdots   \cr
                                            \langle {\cal P}(A), h_{\nu (d)}\rangle _ \beta
                                           \cr} .
\]

To complete the proof we only have to show that $\langle A, h_i\rangle _\beta = 
\langle {\cal P} (A), h_i \rangle _\beta$ for every observable $A$ on $M$.
But, the monomials $h_{o({\bf n})} = H_1^{n_1} H_2^{n_2} \ldots H_k^{n_k}$
are constant on every level set ${\cal L}_\alpha$,  and so
\begin{eqnarray}
\langle A, h_i \rangle _\beta & = & \int _M A(m) h_i(m)  \ e^{-\beta H(m)} \ \dd m  \nonumber\\
& = & \int  _D  h_i^{\cal H}(\alpha) \ e^{-\beta H^{\cal H}(\alpha)}
\Bigl( \int _{{\cal L}_\alpha} A(m_\alpha) \ \dd m_\alpha \Bigr) \ \dd \alpha   \\
& =  &\int _D A^{\cal H} (\alpha)  h_i^{\cal H}(\alpha) \ e^{-\beta \alpha _1} \ \dd \alpha  \nonumber \\
 & = &  \langle A^{\cal H}, h_i ^{\cal H} \rangle _D \nonumber .
\end{eqnarray}
From (\ref{dbeta}) and from the fact that ${\cal P}(h_i) = h_i$,  we now finally get
\[ \fl
\langle A^{\cal H}, h_i^{\cal H} \rangle _D = \langle {\cal P}(A), {\cal P}(h_i) \rangle _\beta =
\langle {\cal P}(A), h_i\rangle _\beta = 
\langle {\cal P}(A), H_1^{m_1}H_2^{m_2}\ldots H_k^{m_k} \rangle _\beta,
\]
 which concludes the proof.
 
 \epf
 
We shall now clarify the question, when  the inequalities in (\ref{new1}) and in 
(\ref{infinite}) are  saturated. The answer is given by  the theorem \ref{corsat} bellow, which is 
more or less an immediate corollary
of theorem \ref{theorem2}.

\begin{thm} {\bf (A) }
Suppose an observable $A$  on $M$ satisfies the two  equivalent conditions:
\begin{enumerate}
\item
The  function $A ^{\cal H} \colon D \to  \rf$ is a polynomial of 
degree $d$ in the variables $(\alpha _1, \ldots , \alpha _k)$.
\item
The observable $A$ be expressible in the form
\begin{equation} \fl
A(m) = \sum _{o({\bf n}) = 0}^{\nu (d)} c_{o({\bf n})}(m) \ H_1^{n_1}(m) \cdots H_k^{n_k}(m) ,
\label{poly}
\end{equation}
where    
\begin{equation} \fl
c_{o({\bf n})}^{\cal H} (\alpha) = \frac{1}{{\rm Vol}({\cal L}_\alpha)} \int _{{\cal L}_\alpha} 
c_{o({\bf n})}(m_\alpha) \ \dd m_\alpha
\equiv const. \quad for \ every  \ {\bf n}.
\label{hconst}
\end{equation}
\end{enumerate}
 Then 
\begin{equation} \fl 
{\cal C}(A)   = \sum _{o({\bf l}), o({\bf n }) = 0}^{\nu (d)}
\langle A  ,   H_1^{l_1} H_2^{l_2} \ldots H_k^{l_k} \rangle _\beta \ 
({\bf H} ^{-1})_{o({\bf l}), o({\bf n})} \ 
\langle A  ,   H_1^{n_1} H_2^{n_2} \ldots H_k^{n_k} \rangle _\beta.
 \label{new1eq}
 \end{equation}
{\bf (B)}  If  $A^{\cal H} \colon D \to \rf$ is an {\rm analytic} function, or alternatively, if
$A$ is expressible in the form
\[ \fl
A(m) = \sum _{o({\bf n}) = 0}^{\infty} c_{o({\bf n})}(m) \ H_1^{n_1}(m) \cdots H_k^{n_k}(m) ,
\]
where the coefficients $c_{o({\bf n})}$ again satisfy the condition  (\ref{hconst}), then
we have the equality
\begin{equation} \fl
{\cal C}(A)   = \sum _{o({\bf l}), o({\bf n }) = 0}^{\infty}
\langle A  ,   H_1^{l_1} H_2^{l_2} \ldots H_k^{l_k} \rangle _\beta \ 
({\bf H} ^{-1})_{o({\bf l}), o({\bf n})} \ 
\langle A  ,   H_1^{n_1} H_2^{n_2} \ldots H_k^{n_k} \rangle _\beta.
 \label{new1eqInf}
 \end{equation}
\label{corsat}
\end{thm}

\noindent {\bf Proof.}
First we check the fact that the conditions {\it 1.} and {\it 2.} of the theorem are 
indeed equivalent. Let the observable $A$ be expressible as
\[
A(m) = \sum _{o({\bf n}) = 0 }^{\nu (d')}  d_{o({\bf n})}(m) H_1^{n_1}(m) \cdots H_k^{n_k}(m).
\]
For every $m$, such that  ${\cal H}(m) = (\alpha_1, \ldots , \alpha_k)$,
we  have $H_i(m) \equiv \alpha _i$ 
on ${\cal L}_\alpha$.  Therefore,
\[
A^{\cal H}(\alpha) = \sum _{ o({\bf n}) = 0}^{\nu (d')} d^{\cal H}_{o({\bf n})}(\alpha) \ 
\alpha _1^{n_1} \cdots \alpha _k^{n_k}.
\]
This function is a polynomial precisely when all $d^{\cal H}_{o({\bf n})}(\alpha)$
are polynomials. In such cases the above function can be rewritten in the form
\[
A^{\cal H}(\alpha) = \sum _{ o({\bf n}) = 0}^{\nu (d)} \lel{c}_{o({\bf n})} \ 
\alpha _1^{n_1} \cdots \alpha _k^{n_k},
\]
where the constants $\lel{c}_{o(\bf{n})}$ 
are coefficients of the  polynomials $d^{\cal H}_{o({\bf n})}(\alpha)$. 
Clearly, every $d_{o({\bf n})}$ is of the form 
$d_{o({\bf n})}(m)  = \sum c_j(m) H_1^{r_1}(m) \cdots H_k^{r_k}(m)$ for
some choice of the multi-indices $(r_1, \ldots r_k)$,  and for every
$c_j(m)$ we have $c_j^{\cal H} = \lel{c}_j$, which proves the equivalence
of {\it 1.} and {\it 2.}

To establish the validity of {\bf (A)}, let 
\[
A^{\cal H}(\alpha_1, \ldots, \alpha _k) = \sum _{o({\bf n }) = 0}^{\nu (d)}
\lel{c}_{o(\bf n)} \ \alpha _1^{n_1} \alpha _2^{n_2} \ldots \alpha _k^{n_k}
\]
be a polynomial.
Then the pull-back of $A^{\cal H}$ on $M$ is given by
\[ \fl
 ({\cal L}(A^{\cal H}))(m) = ({\cal P}(A))(m)  = \sum _{o({\bf n }) = 0}^{\nu (d)}
\Bigl(c_{o(\bf n)} \ H _1^{n_1} H _2^{n_2} \ldots H _k^{n_k}\Bigr)(m)  =
\sum _{o({\bf n}) = 0}^{\nu(d)} c_{o(\bf n)}(m) h_{o({\bf n})}(m) .
\]
According to (\ref{decomp}) this means
that ${\cal P}(A)$ lies in the subspace $V_d$ of the space $L_\beta^2(M)$,
and thus
\[
\pi({\cal P}(A)) = {\cal P}(A).
\]
As seen above, we then have
\begin{eqnarray*}  
{\cal C}(A) & = & \|{\cal P}(A) \|^2_\beta = \| \pi ({\cal P}(A))\|^2_\beta \\
 & = & \sum _{o({\bf l}), o({\bf n }) = 0}^{\nu (d)}
\langle A  ,   H_1^{l_1} H_2^{l_2} \ldots H_k^{l_k} \rangle _\beta \ 
({\bf H} ^{-1})_{o({\bf l}), o({\bf n})} \ 
\langle A  ,   H_1^{n_1} H_2^{n_2} \ldots H_k^{n_k} \rangle _\beta.
\end{eqnarray*}
Part {\bf (B)} of the theorem is an easy consequence of part {\bf (A)}.

\epf

Bellow we collect some remarks and comments that illustrate the condition (\ref{hconst}).


\noindent {\bf Remarks} {\it  
1. \ The coefficients $c_{o({\bf n})}(m)$ in the expression (\ref{poly}) are, of course, in
general not constants. They can be rather arbitrary functions that change in the fibre direction
of the moment map ${\cal H} \colon M \to D$, as well as in the directions transversal to the
fibres ${\cal L}_\alpha$.

\noindent 2. Let $\alpha \in D \subset \rf ^k$ be a regular value of the moment map 
${\cal H} \colon M \to D$.  Then there exists a neighbourhood $\alpha \in U \subset D$
such that  the open subset ${\cal H}^{-1}(U) \subset M$  is diffeomorphic
to $U \times {\cal L}$ and ${\cal L} $ is diffeomorphic to ${\cal L}_\alpha$. Open sets
$W \subset {\cal H}^{-1}(U) \cong U \times {\cal L}$ can be coordinatized as
$W = \{(\alpha _1, \ldots , \alpha _k, \beta _1, \ldots , \beta _l)\}$,
where $(\alpha_1, \ldots , \alpha _k)$ are coordinates on $U \subset D$
and $(\beta_1, \ldots \beta_l)$ some local coordinates on a patch of  ${\cal L}$. 
Then the conditions
\begin {equation}
\frac{\partial}{\partial \alpha _i}\  c_{o({\bf n})} \equiv 0, \quad \quad i = 1, \ldots , k, 
\ \ o({\bf n}) = 0, \ldots , \nu (d)
\label{coords}
\end{equation}
are sufficient for (\ref{hconst}) and hence for (\ref{new1eq}). Clearly, the above
conditions are not necessary for (\ref{hconst}). The functions 
$c_{o({\bf n})}$ are allowed to vary in the $\alpha$-directions, only
their averages over ${\cal L}_\alpha$ have to be constant.

\noindent 3.  Let $M$ be an almost K\" ahler manifold. This means that it is equipped
with a metric $g(-, -)$ compatible with the symplectic form in the sense that 
there exists an almost complex structure $J$ on $M$ such that
\[
g_m(X_m, Y_m) = \omega_m(X_m, J_m(Y_m)), \quad \quad X_m, Y_m \in T_mM.
\]
The most common examples of such manifolds are cotangent bundles $M= T^*N$
over Riemannian manifolds $N$. The metric $g_N$ on $N$ is extended in the natural
way onto the metric on the tangent bundle $TN$,  and the symplectic form $\omega ^T$
on $TN$ is given as the exterior derivative $ \omega ^T= \dd \theta$ of the
tautological 1-form
 \[
 \theta _m (v_m) = (g_N)_{\pi (m)} \Bigl( (D_m\pi)(v_m), m\Bigr), \quad \quad
 v_m \in T_m(TN).
 \]
 Here $\pi \colon TN \to N$ is the natural projection. (Note that $m$ is
  a tangent vector, $m \in T_{\pi(m)}N$.) One can use the metric again to pull
  the form $\omega^T$ back to the cotangent bundle $T^*N$.
Then we have the well defined gradient vector fields $\nabla H_i$ associated to functions 
$H_i$, which are $g$-orthogonal to the Hamiltonian vector fields $X_H$.
Moreover, the gradients $\nabla H_i$ are orthogonal to the entire fibres
${\cal L}_\alpha$. Therefore, the conditions
\[
(\nabla H_i)(c_{o({\bf n})}) \equiv 0, \quad i = 1, \ldots , k, \ \  o({\bf n}) = 0, \ldots , \nu (d)
\]
are equivalent to the conditions (\ref{coords}).  Similarly as the conditions
(\ref{coords}),  they are sufficient but
not necessary for (\ref{hconst}).

\noindent 4. \  The condition that $A^{\cal H}$ is a polynomial is also equivalent to the condition
\[
if \ |{\bf n}| = n_1 + n_2 + \ldots +  n_k > d, \ \quad  then \quad
\langle A, H_1^{n_1} H_2^{n_2} \ldots H_k^{n_k} \rangle _\beta = 0.
\]
This follows immediately from the expansion (\ref{decomp}) and from the equivalence of the
conditions 1. and  2. of  theorem \ref{corsat}.  In practice, the problem with the 
above condition is that infinitely many integrals have to be checked.
But on the other hand, in some cases, the integrals of the functions $A \cdot  H_1^{n_1} \cdots  H_k^{n_k}$
might be rather easily computable in some concrete cases.}


\section{The general case}

The ergodically regular systems treated above are very special. The invariant subspaces
on which a Hamiltonian system is ergodic can in general be rather wild, and even within
a single system they can be of very different types, e.g.   full level 
sets ${\cal L}_\alpha$, subsets of ${\cal L}_\alpha$ of lower dimensionality, as well as 
subsets of ${\cal L}_\alpha$ having even fractal dimensions. 
One should recall for example 
typical situations of smooth perturbations of integrable systems landing in the context of the 
Kolmogorov-Arnold-Moser theory.
In general,  it is impossible to parameterize
invariant ergodic sets by some manageable (say Hausdorff) space over which one could integrate. 
For these reasons we have to modify our approach in order to prove our bounds for a general 
Hamiltonian system.  In particular, for a useful description of the ergodic
decomposition of our system, we shall revert to an inverse limit type construction. This will
enable us to prove the results analogous to those from the previous section, but valid
for the general Hamiltonian systems. The analogue of theorem \ref{theorem1} is the following:

\begin{thm}
Let $(M, \omega, H)$ be an arbitrary Hamiltonian system with a well defined partition
function and let $A \colon M \to \rf$ be an element of $L_\beta^1(M)$.
Then for the orbital average 
\[
\lel{A}(m) = \lim _{ T \to \infty} \frac{1}{T} \int _0 ^T A(\gamma _m(t)) \ \dd t
\]
we have
\[
{\cal C}(A) = \| \lel{A}\|^2_{\beta}.
\]
\label{theorem3}
\end{thm}

Before giving the proof, we shall describe the ergodic decomposition of the 
Hamiltonian system $(M, \omega, H)$ which will be used in the proofs.
In our construction we shall use the invariant
measure on $M$ given by
\[
\mu (N) = \frac{1}{{\cal Z}(\beta)} \ \int _N e^{-\beta H(m)} \ \dd m .
\]
As in the previous section, the crux of the proof will be the replacement of the
temporal averages by the spatial ones in the context without ergodic regularity.
To be able to do this, we have to decompose
the space $M$ into a collection of invariant sets on which our system 
{\em is} ergodic.  We will construct such ergodic decomposition by means 
of successive approximations.
By definition the Hamiltonian system $(M, \omega, H)$ is {\it not} ergodic, if there exists an invariant
measurable set  $N \subset M$  such that
\[
0 < \mu (N) < 1,
\]
where the inequalities have to be strict. Let the first approximation of
our ergodic decomposition be a finite partition ${\cal N}_1$ consisting of {\em invariant sets}
\[
N_1^1, N_2^1, \ldots , N_{k_1}^1  \subset M
\]
with the properties
\[
0 < \mu(N_i^1) < 1, \quad \mu(N_i^1 \cap N_j^1) = 0, \quad \cup _{i =1}^{k_1} N_i^1 = M.
\]
In the next stage we decompose each $N^1_i$ into {\it invariant} measurable sets
$\{N^2_{k_i}, N^2_{k_i + 1},  \ldots , N^2_{l_i}\}$ with analogous properties.
This yields the partition ${\cal N}_2 = \{N^2_1, \ldots N_{k_2}^2\}$ which again satisfies the 
stipulations analogous to those listed above.
We continue the procedure and obtain a sequence $\{ {\cal N}_n\}_{n \in \mathbb{N}}$
of partitions in which every term ${\cal N}_n = \{N_1^n, \ldots , N_{k_n}^n\}$ satisfies 
\[
0 < \mu(N^n_i) < 1, \quad \mu(N_i^n \cap N_j^n) = 0, \quad \cup _{i = 1}^{k_n} = M.
\]
Let now $\{a_n\}_{n \in \mathbb{N}}$ be a sequence of natural numbers such that
for every $i \in \mathbb{N}$ we have
\[
1 \leq a_i \leq k_i.
\]
Then  the set
\[
N_{\{a_n\}} = \lim _{n \to \infty} N^n _{a_n} = \cup _{n = 1}^{\infty} N^n _{a_n}
\]
is measurable and  invariant. Moreover, the system $(M, \omega, H)$ is ergodic
on every $N_{\{a_n\}} $.  


\noindent {\bf Remark} 
{\it The sets $N_{\{a_n\}}$ are subsets of  level sets ${\cal L}_\alpha$, therefore
they are sets of measure  zero. These sets 
can be strange, of noninteger Hausdorff dimensions etc... 
The collection of $N_{\{a_n\}}$  is parameterized by some subset of the set of  sequences
$\{a_n\}_{n \in \mathbb{N}}$ of the form described above.  In order to simplify the 
parameter set, one could make the above  construction in a ``binary manner" by
decomposing each set $N^n_i$ into only two invariant sets $N^{n + 1}_0$ and
$N^{n + 1}_1$. The sequences $\{a_n\}_{n \in \mathbb{N}}$ would then be
the maps 
\[
\{a_n\}_{n \in \mathbb{N}} : \mathbb{N} \longrightarrow \{ 0, 1\}.
\]}


Now we construct the sequence $\{A_n\}_{n \in \mathbb{N}}$ of measurable functions
\[
A_n : M \longrightarrow \rf
\]
by the rule
\begin{equation}
 \quad
if  \quad m \in N_i^n \quad then \quad A_n(m) = 
\frac{1}{{\rm Vol}(N_i^n)} \ \int _{N_i^n} A(h) \ \dd h
= C_{n, i}.
\label{anav}
\end{equation}

\noindent {\bf Proof of Theorem \ref{theorem3}:} Our strategy here will be 
analogous to the one used in the proof of  theorem {\ref{theorem1}.
If $A \colon M \to \rf$ is an element
of $L_\beta^1(M)$, then 
by Birkhoff's theorem (see,  e. g. \cite{mane}) the function $\lel{A}$ is also measurable and $ \lel{A} \in L^1_\beta (M)$.
We clearly have
\[
\lim _{T \to \infty} \frac{1}{T} \int _0 ^T \ \lel{A} (\gamma _m (t)) \ \dd t =
\lel{A} (m),    
\]
which means that the operator
\[
\begin{array}{cccl}
{\cal A} : & L^2 _\beta (M) & \longrightarrow  & L^2 _\beta (M) \\
 & & & \\
  & A &  \longmapsto  & \lel{A}
\end{array}
\]
is a projector.  The image of ${\cal A}$ are the those functions in $L^2_\beta(M)$ which 
are constant on every $H$-orbit. This is a closed subspace in $L^2_\beta(M)$, therefore 
${\cal A}$ is a continuous operator. Moreover, it is also an {\it orthogonal} projector. To see this,
let $B \in {\rm ker}{\cal A}$,  and  $\lel{A}  = {\cal A}(A)$ be arbitrary elements in the 
kernel and in the image of ${\cal A}$ respectively. We claim that
\begin{equation}
\langle B, \lel{A} \rangle _\beta =  \int _M B(m) \ \lel{A}(m) \ \dd \mu = 0.
\label{orthog}
\end{equation}
Let $\{A_n\}_{n \in \mathbb{N}}$ be the sequence of functions approximating $A$ 
as described above. By construction we have: For every $m \in M$ there exists
a unique sequence $\{a_n\}_{n \in \mathbb{N}} $ such that
$m \in N_{\{a_n\}} = \lim _{n \to \infty}N_{a_n}^n$.  Therefore,
\[
\lel{A}(m) = \lim _{n \to \infty} A_n(m),
\]
since our system is ergodic on every $N_{\{a_n\}}$.
This implies 
\[
B(m) \lel{A}(m) = \lim _{n \to \infty} B(m) A_n(m).
\]
Functions $A$ and $B$ are elements of $L_\beta^1$. 
Let now $Sup$ be the essential supremum of the orbital average $\lel{A}$ on $M$. 
The measure $\mu$ on $M$ is finite, so a  function $Sup$ 
is an  element of $L_\beta^1(M)$.  We have the inequality
\[
B(m) A_n(m) \leq B(m) \ Sup
\]
which holds for almost every
$m \in M$. Therefore, by the Lebesgue dominated convergence theorem
we have
\begin{equation}
\int _M B(m) \ \lel{A}(m) \ \dd \mu = \lim _{n \to \infty} \int _M B(m) \ A_n(m) \ \dd \mu.
\label{lebesgue}
\end{equation}
Since $A_n$ takes the  constant value $C_{i, n}$ on every $N_i^n$ for $i = 1, \ldots, k_n$, 
and since all the sets $N_i^n$  as well as the measure $\mu$ are $H$-invariant,
the Liouville theorem gives
\begin{eqnarray}
\int _M B(m) \ A_n(m) \ \dd \mu & = & \sum _{ i = 1}^{k_n} C_{i, n} \int _{N_i^n} B(m) \ \dd \mu 
\nonumber \\
& = &  \sum _{ i = 1}^{k_n} C_{i, n} \int _{N_i^n} B(\gamma _m(t)) \ \dd \mu
\end{eqnarray}
for every $t \in \rf$. Since ${\cal A}(B) = \lel{B}  = 0$,  we get
\[
\int _M B(m) \ A_n(m) \ \dd \mu = \sum _{ i = 1}^{k_n} C_{i, n} 
\lim _{T \to \infty} \frac{1}{T} \int _0^T\Bigl( \int _{N_i^n} B(\gamma _m(t)) \ \dd \mu \Bigr)
\ \dd t = 0,
\]
which together with (\ref{lebesgue}) proves (\ref{orthog}).

Now, every orthogonal projection is a symmetric operator, therefore
\begin{eqnarray*}
{\cal C}(A) & = & \lim _{T \to \infty} \frac{1}{T} \Bigl( \int _M A(m) \ A(\gamma _m(t)) \ \dd \mu \Bigr)
\ \dd t  \nonumber \\
& = &\langle A, {\cal A} (A) \rangle _\beta  = 
\langle A , {\cal A}^2 (A) \rangle _\beta \\
& = & \langle {\cal A}(A)  , {\cal A}(A) \rangle _\beta  \nonumber \\
& = & \| \lel{A} \| ^2_{\beta} \nonumber ,
\end{eqnarray*}
which completes the proof.

\epf

In most cases the orbital average $\lel{A}$ is impossible to calculate. Therefore, theorem \ref{theorem3}
almost never provides a good estimate for ${\cal C}(A)$, apart from the fact that it ensures
${\cal C}(A)$ to  be nonnegative.  Therefore, as we have done in the 
ergodically regular case, we shall construct more easily calculable
estimates which will use the specific information about the Hamiltonian system in question. 
In the general case this information has two sources, the conserved quantities and some
ergodic decomposition $\{ {\cal N }_n \}_{n \in \mathbb{N}}$ of the form described above.

In the proof of theorem \ref{theorem3} we used the sequence $\{A_n\}_{n \in \mathbb{N}}$
of functions approximating the observable $A$.  We shall now  modify the approximating
sequence, so that it will take into account the conserved quantities $H_1, \ldots , H_k$
of the system as well as the ergodic decomposition. We shall mimic the approach 
from the previous section, but will replace the averaging over the level sets
${\cal L}_\alpha$ with the averaging over the intersections
\[
N_{\alpha, i}^n = {\cal L}_\alpha \cap N_i^n .
\]
Consider the functions $B_{n, i} \colon M \to \rf$ for $n\in \mathbb{N}$ and $i \in \{1, \ldots k_n\}$
associated to the ergodic decomposition $\{{\cal N}_n\}_{n \in \mathbb{N}}$,  and given by the rule

\begin{equation}
B_{n, i}(m) = \left\{ 
\begin{array}{ccc}
  \frac{1}{{\rm Vol}(N^n_{{\cal H}(m), i})} \int _{N^n_{{\cal H}(m), i}}  \! A(m_\alpha) 
\ \dd m_\alpha &  \ ; \  \ & m \in N_i^n  \\
& & \\
  0  & \  ; \  \ &   otherwise 
\end{array} . 
\right .
\label{defBn}
\end{equation}
Let  now the sequence $\{C_n\}_{n \in \mathbb{N}}$ of functions 
\[
C_n : M \longrightarrow \rf
\]
be given by
\begin{equation}
C_n(m) = \sum _{i = 1}^{k_n} B_{n, i}(m).
\label{defCn}
\end{equation}
The functions $C_n$ are approximations of the orbital average of $A$. In addition, their behaviour
is similar to that of the projection ${\cal P}(A)$ in the ergodically regular case.
More precisely, the functions $C_n$ are refinements of ${\cal P}(A)$ which take
into account  the partition  ${\cal N}_n$.

Recall that
\[
({\cal P}(A))(m) = \frac{1}{{\rm Vol}({\cal L}_{{\cal H}(m)})} \int _{{\cal L}_{{\cal  H}(m)}}
A(m_\alpha) 
\ \dd m_\alpha, \quad for \ m \in {\cal L}_\alpha.
\]
Since ${\cal L}_\alpha = \cup _{i = 1}^{k_n} N^n_{\alpha, i}$, for every  $m \in {\cal L}_\alpha$, we have
\begin{eqnarray}
({\cal P}(A))(m) & = & \frac{1}{{\rm Vol}({\cal L}_\alpha)}  \sum _{i = 1}^{k_n}
\int _{N^n_{\alpha, i}} A(m_\alpha)  \ \dd m_\alpha   \nonumber \\
& = & \sum _{i =1}^{k_n}  C_n(m_i) \ \frac{{\rm Vol}(N^n_{\alpha, i})}{{\rm Vol}({\cal L}_\alpha)}
\quad \quad m_i \in N_{{\cal H}(m), i}^n \ arbitrary . 
\label{avsum}
\end{eqnarray}

\begin{lema}
The sequence of functions $\{C_n\}_{n \in \mathbb{N}}$ has the following properties.
\begin{enumerate}
\item For every $n \in \mathbb{N}$ we have
\[
{\cal C}(C_n) = \|C_n\|^2_\beta \geq \| {\cal P}(A) \|^2_\beta = \| A^{\cal H}\|^2_D.
\]
\item The sequence $\{\| C_n\|^2_\beta \}_{n \in \mathbb{N}}$ is non-decreasing.
\item The sequence $\{\| C_n\|^2_\beta \}_{n \in \mathbb{N}}$ is convergent, and
\[
{\cal C}(A)  = \lim _{n \to \infty} {\cal C}(C_n) = \lim _{ n \to \infty} \| C_n\|^2_\beta. 
\]
\end{enumerate}
\label{lemma}
\end{lema}
{\bf Proof.}
{\it Ad 1.}
The definition of the time average of the correlation function, and the fact that 
the functions $C_n$ are constant on the $H$-invariant sets $N_{\alpha, i}^n$ give
\begin{eqnarray*}
{\cal C}(C_n) & = & \int _M C_n(m) \Bigl(\lim_{T \to \infty}\frac{1}{T} 
\int _0^T C_n(\gamma _m(t)) \ \dd t \Bigr) \ \dd \mu \\
& = &\int _D  \dd \alpha \sum _{ i = 1}^{k_n} \int _{N_{\alpha, i}}
B_{n, i}(m_\alpha)  \Bigl(\lim_{T \to \infty}\frac{1}{T} 
\int _0^T B_{n, i}(\gamma _{m_{\alpha}}(t)) \ \dd t \Bigr) \ \dd m_\alpha \\
& = &\int _D \dd \alpha \sum _{i = 1}^{k_n} \int _{N_{\alpha, i}^n} 
B^2_{n, i}(m_\alpha)  \ \dd m_\alpha .
\end{eqnarray*}
On the other hand,  we have
\begin{eqnarray}
\|C_n\|^2_\beta  & = & \int _M C_n(m) ^2 \ \dd \mu  \nonumber \\
& = &
\int _D \ \dd \alpha \sum _{i =1}^{k_n} \int _{N_{\alpha, i}^n} B^2_{n, i}(m_\alpha) \ \dd m_\alpha,
\label{norm}
\end{eqnarray}
which proves the first equality in {\it 1.}

To prove the inequality $\|C_n\|^2_\beta \geq \|{\cal P}(A)\|_\beta^2$ we use (\ref{avsum})
and the usual procedure for calculating variances. This gives
\begin{eqnarray*}
0  & \leq & \sum _{ i = 1}^{k_n} \Bigl(   C_n(m_i)  - {\cal P}(A) (m)  \Bigr) ^2 \ 
\frac{{\rm Vol}(N_{\alpha, i}^n)}{{\rm Vol}({\cal L}_\alpha)}  \\
& = & \sum _{i = 1}^{k_n} C_{n}^2(m_i) \frac{{\rm Vol}(N_{\alpha, i}^n)}{{\rm Vol}({\cal L}_\alpha)}
 - {\cal P}(A)^2(m) 
 \quad \quad m_i \in N_{\alpha, i}^n \  \ arbitrary
\end{eqnarray*}
for every $ m \in {\cal L}_\alpha$.   Integrating the above inequality over $M$ with respect
to the measure $\dd \mu$ gives
\[
\int _M \Bigl( \sum _{i = 1}^{k_n} C_{n}^2(m_i) 
\frac{{\rm Vol}(N_{\alpha, i}^n)}{{\rm Vol}({\cal L}_\alpha)} \Bigr) \ \dd \mu  \ \geq \ 
\int _M {\cal P}(A)^2(m)  \ \dd \mu = \| {\cal P}(A)\|^2_\beta .
\]
Observing that the expression 
$\sum _{i = 1}^{k_n} C_{n}^2(m_i) 
\frac{{\rm Vol}(N_{\alpha, i}^n)}{{\rm Vol}({\cal L}_\alpha)}$
is an average of a function over ${\cal L}_\alpha$ and is therefore constant
on ${\cal L}_\alpha$ yields
\begin{eqnarray*}
\int _M \Bigl( \sum _{i = 1}^{k_n} C_{n}^2(m_i) 
\frac{{\rm Vol}(N_{\alpha, i}^n)}{{\rm Vol}({\cal L}_\alpha)} \Bigr) \ \dd \mu & = &
\int _D \dd \alpha \int _{{\cal L}_\alpha} \sum _{i = 1}^{k_n} C_{n}^2(m_i) 
\frac{{\rm Vol}(N_{\alpha, i}^n)}{{\rm Vol}({\cal L}_\alpha)} \ \dd m_\alpha \\
& = & \int _D  \Bigl( \sum _{i = 1}^{k_n} C_{n}^2(m_i) 
{\rm Vol}(N_{\alpha, i}^n) \Bigr)  \ \dd \alpha .
\end{eqnarray*}
Since $C_n(m_i) = B_{n, i}(m)$ for suitable pairs  of $m$ and $m_i$, we see from
(\ref{norm}) that the above expression is indeed equal to the 
norm $\|C_n\|^2_\beta$, which concludes the proof of the point {\it 1}.

\noindent {\it Ad 2.}
The proof that for every $n\in \mathbb{N}$ we have $\|C_{n + 1}\|_\beta^2 \geq \| C_n\|^2_\beta$,
is essentially the same as the proof of the inequality
$\|C_n \|^2_\beta \geq \|{\cal P}(A)\|^2_\beta$ just given.
Above we partitioned the phase space into the disjoint union $\cup _{i = 1}^{k_n} N_i^n = M$
of $H$-invariant subsets with positive measures.  To prove {\it 2.} we have to partition
every $N_i^n$ further into the union $\cup _{j = l_i}^{k_i} N_j^{(n + 1)} = N_i^n$.
The actual calculations here are then precisely the same as in {\it 1.},  modulo slightly different
notation.

\noindent {\it Ad 3.} First we observe that for every $m \in M$ we have
\[
\lim _{ n \to \infty} C_n(m) = \lel{A}(m) .
\]
Indeed, for every $m$ there exists a unique $\alpha$ such that $m \in {\cal L}_\alpha$, 
and a unique sequence $\{a_n\}_{n \in \mathbb{N}}$
such that $m \in N_{\{a_n\}} = \cap _{n = 1}^{\infty} N_{a_n}^n$.  Since our
system is ergodic on the limit set $N_{\{a_n\}}^n = \lim _{n \to \infty} N_{a_n}^n$, 
and since the orbit $\gamma _m(t)$ is contained in ${\cal L}_\alpha$, we have
\begin{eqnarray*}
\lim _{n \to \infty} C_n(m) & = & \lim _{n \to \infty} \frac{1}{{\rm Vol}(N^n_{\alpha, a_n})}
\int _{N_{\alpha,  a_n }^n} A(m_\alpha) \ \dd m_\alpha \\
& = & \lim _{T \to \infty}\frac{1}{T} \int _0^T A(\gamma _m(t)) \ \dd t  \\
& = & \lel{A}(m).
\end{eqnarray*}
Let again $Sup$ be the essential supremum of the orbital average $\lel{A}$ on $M$. Then
clearly 
\[
C_n(m) \leq Sup; \quad a. \ e. \ on \ M ,
\]
since the values $C(m)$ are averages taken over larger sets than those that yield $\lel{A}(m)$.
As we already mentioned, $Sup$ is an element of $L_\beta^1(M)$.
Thus, the sequence $\{ C_n^2\}_{n \in \mathbb{N}}$ together with its
point-wise limit $\widetilde{A}^2$ satisfies the conditions of the Lebesgue dominated 
convergence theorem. Therefore,  we have
\[
\lim _{n \to \infty} \| C_n\|^2_\beta = \lim _{n \to \infty} \int _M C_n^2(m) \ \dd \mu
= \int _M \lel{A}^2(m) \ \dd \mu = \|\lel{A}\|^2_\beta.
\]
This, together with the result of theorem \ref{theorem3}, concludes the proof.

\epf

An immediate consequence of lemma \ref{lemma} is the fact that all the inequalities
proved in the previous section for the ergodically regular systems hold for all
Hamiltonian systems without restrictions. In particular, we have
\begin{cor}
Let $(M, \omega, H)$ be an arbitrary Hamiltonian system with $k$ conserved
quantities $H = H_1, H_2, \ldots H_k$. Then for every observable 
$A \colon M \to \rf$ and every $d \in \mathbb{N}$ we have:
\begin{equation} \fl
{\cal C}(A)   \geq \sum _{o({\bf l}), o({\bf n }) = 0}^{\nu (d)}
\langle A  ,   H_1^{l_1} H_2^{l_2} \ldots H_k^{l_k} \rangle _\beta \ 
({\bf H} ^{-1})_{o({\bf l}), o({\bf n})} \ 
\langle A  ,   H_1^{n_1} H_2^{n_2} \ldots H_k^{n_k} \rangle _\beta.
 \label{new1eqGen}
 \end{equation}
 The above inequality also holds for $d = \infty$.
\label{corlast}
\end{cor}
{\bf Proof.}
We have just proved that for a general system  we have ${\cal C}(A) \geq {\cal P}(A)$. Proposition
then follows from theorem \ref{theorem2}.

\epf

We notice that the estimate in the above proposition does not reflect in any way how far
our system is from being ergodic on the level sets ${\cal L}_\alpha$. The non-ergodicity is
reflected in the ergodic decomposition $\{ {\cal N}_n\}_{n \in \mathbb{N}}$,  and the information
given by $\{ {\cal N}_n\}_{n \in \mathbb{N}}$ can be used to improve the bound
(\ref{new1eqGen}). In our setup it is quite easy to plug the decomposition ${\cal N}_n$ into
 (\ref{new1eqGen}).  Recall the definitions (\ref{defBn}) and (\ref{defCn}).
Since  for $i \ne j$  the supports $N_i^n$ and $N_j^n$ of 
the functions $B_{n, i}$  and $B_{n, j}$ are disjoint,  we have
\begin{equation}
\langle B_{n, i}, B_{n, j} \rangle _\beta = 0.
\label{orto}
\end{equation}
The function $C_n = \sum _{i = 1}^{k_n} B_{n, i}$ is a sum of orthogonal vectors, 
therefore, 
\begin{equation}
\| C_n\|^2 _\beta = \sum _{ i = 1}^{k_n} \| B_{n, i}\|^2_\beta.
\label{pyth}
\end{equation}
As in the previous section, we can project orthogonally  the function $C_n \in L_\beta^2(M)$
on the subspace $V_d \subset L_\beta^2(M)$ spanned by the monomials 
$H_1^{n_1} H_2^{n_2} \ldots H_k^{n_k}$ of degree $d$ or less. We have 
proved  that the vector $C_n$ is longer than ${\cal P}(A)$, therefore for 
the projections by $\pi \colon L_d^2(M) \to V_d$ we have
\begin{equation}
\pi(C_n) \geq \pi ({\cal P}(A)).
\label{leqproj}
\end{equation}
Taking into account (\ref{orto}) and (\ref{pyth}),  and lemma \ref{lemma},
we get the following proposition.
\begin{prop}
For every $n \in \mathbb{N}$ we have
\[ \fl
{\cal C}(A) \geq  \sum _{ i = 1}^{k_n} \Bigl(
\sum_{o({\bf l}), o({\bf j}) = 0}^{\nu (d)} \langle B_{n, i}, H_1^{l_1} H_2^{l_2} \ldots H_k^{l_k} 
 \rangle _\beta \cdot ({\bf H}^{-1})_{o({\bf l}), o({\bf j})} \cdot
 \langle B_{n, i}, H_1^{j_1} H_2^{j_2} \ldots H_k^{j_k} \rangle _\beta \Bigr).
\]
The inequality also holds for $d = \infty$.
\label{prop4}
\end{prop}
\epf

\noindent The quality of the above estimate increases with increasing $n$ and $d$,
and for every $n \geq 2$ the above estimate  is better that the estimate in  theorem
\ref{theorem2}.

\section*{Acknowledgements}

The authors would like to express their gratitude to Mirko Degli Esposti,
Andreas Knauf and Oliver Dragi\v cevi\' c for their valuable comments
and suggestions. The work has been supported by the programmes P1-0291, P1-0044, and the grant J1-7347 
of the Slovenian Research Agency (ARRS).


\section*{References}

\begin{thebibliography}{1}

\bibitem{mazur} Mazur P 1969 Physica {\bf 43} 533 

\bibitem{suzuki} Suzuki M 1971 Physica {\bf 51} 277

\bibitem{kubo} Kubo R 1957 J. Phys. Soc. Japan {\bf 12} 570

\bibitem{visscher} For a very clear account on Green-Kubo theory in the realm of classical mechanics see: 
Visscher W M 1974 Phys. Rev. A {\bf 10} 2461

\bibitem{khinchine} Khinchine A I 1934 Ann. Math. {\bf 109} 604

\bibitem{zotos} Zotos X, Naef F and Prelov\v sek P 1997 Phys. Rev. B {\bf 55} 11029

\bibitem{zotosreview} Zotos X and Prelov\v sek P, review article published in
``Strong Interactions in Low Dimensions", series ``Physics and Chemistry of Materials with Low Dimensional Structures", eds. D. Baeriswyl and L. Degiorgi, Kluwer Academic Publishers (2004)

\bibitem{prosen98} Prosen T 1998 J. Phys. A: Math. Gen. {\bf 31} L645 

\bibitem{minlos} Minlos R A 2002 {\it Introduction to Mathematical Statistical Physics} (In Russian)
Moscow: MCNMO)

\bibitem{mane} Mane R 1987 {\it Ergodic Theory and Differentiable Dynamics} 
(Berlin: Springer-Verlag)

 \end{thebibliography}
\end{document}